%% file: GrowingSpines.tex
\documentclass{amsart}
\usepackage{graphicx} 

\title{Growing spines ad infinitum}

\author[B.~Boissonneau]{Blaise Boissonneau}
\address{Blaise Boissonneau \\ Heinrich Heine University Düsseldorf, Faculty of Mathematics and Natural Sciences, Universitätsstr.~1, 40225 Düsseldorf, Germany. \url{blaise.boissonneau@hhu.de}}

\author[A.~De Mase]{Anna De Mase}
\address{Anna De Mase \\ Università degli Studi della Campania ``Luigi Vanvitelli'', Department of Mathematics and Physics, viale Lincoln 5, 81100 Caserta, Italy. \url{anna.demase@unicampania.it}}
\thanks{ADM was partially supported by GNSAGA (INdAM, Italian National Institute of High Mathematics) - Project ``Model theory and applications to algebra and combinatorics", and PRIN 2022 ``Models, sets and classification".}

\author[F.~Jahnke]{Franziska Jahnke}
\address{Franziska Jahnke \\ Institute for Logic, Language and Computation, University of Amsterdam, 
Science Park 107,
1098 XG Amsterdam, The Netherlands 
and Institute for Mathematical Logic and Foundations, Department of Mathematics and Computer Science,
University of M\"unster,
Einsteinstraße 62,
48149 M\"unster, Germany. \url{franziska.jahnke@uni-muenster.de}}
\thanks{FJ was supported by the Deutsche Forschungsgemeinschaft (DFG, German Research Foundation) under Germany's Excellence Strategy EXC 2044-390685587, Mathematics M\"unster: Dynamics-Geometry-Structure, as well as by a Fellowship from the Daimler and Benz Foundation.}

\author[P.~Touchard]{Pierre Touchard}

\address{Pierre Touchard \\ Institut f\"ur Algebra, Technische Universit\"at Dresden, 01062 Dresden, Germany and KU Leuven, Department of Mathematics, B-3001 Leuven, Belgium. \url{pierre.touchard@tudresden.de}}

\thanks{PT was partially supported by KU Leuven IF C16/23/010.}

\date{\today}
\input{preamble.tex}

\begin{document}

\begin{abstract}
We show that every non-trivial ordered abelian group $G$ is augmentable 
by infinite elements, i.e., we have $G \substruct H \oplus G$ for some non-trivial 
ordered abelian group $H$. As an application, we show that when $k$ is a
field of characteristic $0$, then $k$ is not $t$-henselian if and only if all henselian
valuations with residue field $k$ are ($\emptyset$-)definable.
\end{abstract}

\maketitle
\tableofcontents

\setcounter{section}{-1}
\section{A note on a former result}

This paper was written without knowledge of the paper ``Inclusions et produits de groupes abéliens ordonnés étudiés au premier ordre'' by Françoise Delon and François Lucas. In this paper from 1989, they prove, among other results, that given $G\substruct G^*$ ordered abelian groups and $H=\langle G\rangle$ the convex hull of $G$, then $G\substruct H\substruct G^*$, that is, our \Cref{lem:ConvxHull} is exactly Proposition 1 from \cite{DelonLucas}.

In order to prove this, they work in the formalism of Schmitt, see \cite{Sch82}. They argue in Example 1 that given coloured linear orders $C_1\substruct C_2$, if $C_3$ is either the initial closure of $C_1$ in $C_2$, the final closure of $C_1$ in $C_2$, or the convex hull of $C_1$ in $C_2$, then $C_1\substruct C_3\substruct C_2$. This is essentially what we prove in \Cref{lem:hulls!}, although we are doing it for multi-orders instead of linear orders.

Their Lemma 1-1 is saying that if $H$ is a convex subgroup of an ordered abelian group $G$, then $Sp_n(H)$ is an initial segment of $Sp_n(G)$, where $Sp_n$ denote the Spectrum of order $n$, a structure coming from Schmitt's formalism, and very similar to the $n$-spines of Cluckers-Halpczok's language. This result is therefore very similar to our \Cref{cor:spines-of-cvx}.

It is then just a matter of piecing it together to obtain Proposition 1 from their paper or \Cref{lem:ConvxHull} of this paper.

Because we didn't know that their result existed, we ended up reproving it ourselves independently, going by a very similar argument but using Cluckers-Halupczok's language instead of Schmitt's.

Our main theorem, \Cref{thm:strongleftaugmentabilityOAG}, is however new, to the extent of our knowledge. Its application to definability of valuations is also a new result.

\section{Introduction}
In this paper, we study the model theory of ordered abelian groups, building on the work of
Schmitt \cite{Sch82} as well as Cluckers and Halupczok \cite{CH11}.
Ordered abelian groups are much-studied objects in model theory since the pioneer work of Robinson and Zakon \cite{RoZa60}. Among important and classical contributions, we may cite \cite{DelonLucas} (on elementary subgroups), \cite{Sch84} (quantifiers elimination relative to spines) and \cite{Kuh90} (ordered and valued divisible groups).
Similar methods were developed in the context of valued vector spaces in \cite{KK95} and applied in \cite{Kuh99}. There are several recent developments
on the classification of ordered abelian groups according to dividing lines (\cite{JSW}, \cite{ACGZ20}, \cite{HH19A}), elimination of imaginaries \cite{Vic22}, and 
pairs of models (e.g. \cite{CHY21}, \cite{HLT23}).
The study of ordered abelian groups is often motivated by that of valued
fields, since by the Ax-Kochen/Ershov principle, understanding value groups is a key step in order to understand valued fields. This connection is explicit in e.g.~\cite{ACGZ20}, \cite{HM21}, and \cite{V23}. Properties of value groups are also often used to define henselian valuations
in the language of rings, see \cite{Hon14} and \cite{KKL}.

By the results of Schmitt \cite{Sch82} and Cluckers-Halupczok \cite{CH11}, many first-order
properties of ordered abelian group can be reduced to corresponding properties of spines, which are chains of uniformly definable convex subgroups.
However, due to the highly technical nature of the framework, 
most of the works cited above restrict to subclasses, e.g., ordered abelian groups with finite spines, and results applying to the class of all ordered abelian groups are rare. 

In this paper, we study augmentability of ordered abelian groups by infinite elements.
We call an
ordered abelian group $G$ \emph{augmentable by infinites} if there exists a non-trivial ordered abelian group $H$ such that $G\substruct H\oplus G$. Our main result is that this holds for 
any non-trivial ordered abelian group:

\begin{theorem}[Cf.~Theorem \ref{thm:strongleftaugmentabilityOAG}]
    All non-trivial ordered abelian groups are augmentable by infinites.
\end{theorem}

The proof of this theorem goes via spines.
We also characterize which ordered abelian groups are augmentable by an ordered abelian group $H$ which is ($p$-)divisible (Corollary \ref{cor_div}).

As an application of our results, we study definable henselian valuations. Here, we show that
if $k$ is a field of characteristic $0$ that is $\mathcal{L}_\textrm{ring}$-elementarily
equivalent to a field admitting a non-trivial henselian valuation (that is, if $k$ is $t$-henselian),
then $k \substruct k((\Gamma))$ for some non-trivial ordered abelian group $\Gamma$ (Proposition \ref{power}).
We apply this to deduce the following:
\begin{theorem}[Cf.~Theorem \ref{def_fieldwise}]
    Let $k$ be a field of characteristic $0$. The following are equivalent:
    \begin{enumerate}
        \item $k$ is not t-henselian,
        \item for every henselian valued field $(K,v)$ with residue field $k$, the valuation ring $\mathcal{O}_v$ is $\mathcal{L}_\textrm{ring}$-definable 
        (possibly using parameters),
        \item for every henselian valued field 
        $(K,v)$ with residue field $k$, the valuation ring $\mathcal{O}_v$ is 
        $\emptyset$-$\mathcal{L}_\textrm{ring}$-definable,
        \item All henselian valuation rings with residue field elementarily equivalent to $k$ are
        uniformly $\emptyset$-$\mathcal{L}_\textrm{ring}$-definable in $\mathcal{L}_\textrm{ring}$.
    \end{enumerate}
\end{theorem}

This partially answers a question by Krapp, Kuhlmann and Link \cite{KKL}, 
who ask for a characterization of the class of fields such that for any
henselian valued field with residue field in said class, 
the valuation is (parameter-freely) definable in the language of rings. 
A similar characterization for existential (respectively universal) definability of the valuation was proven by Anscombe and Fehm in \cite{AF17}.
\smallskip

The paper is organized as follows. In \hyperref[sec:spines]{Section 2}, we give an introduction to spines, as well
as the language $\mathcal{L}_{\text{syn}}$. We study the behaviour of convex subgroups in this language, and discuss when the spines of a convex subgroup $H \subseteq G$ can be embedded into the spines of $G$ (Corollary \ref{cor:cvx-are-substruct}). Spines can be seen as coloured multi-orders, a multisorted generalization of linear orders (see \Cref{def:multiord}).
\hyperref[sec:multiord]{Section 3} hence discusses coloured multi-orders. 
In particular, we show that given an elementary pair $A \substruct B$ of coloured multi-orders, the convex hull of $A$, as well as the left and right closures of $A$, are elementary
substructures of $B$ (Lemma \ref{lem:hulls!}). In the remainder of \hyperref[sec:multiord]{Section 3}, we record some facts
about augmentability of coloured multi-orders for future reference. 
Most notably, in Lemma \ref{lem:aug-mo}, we characterize 
augmentability of a coloured multi-order $A$ on the right (resp.~left):
We show that there is a coloured multi-order $B$ such that $A \substruct A + B$ 
(resp.~$A \substruct B +A$) holds if and only if some component of $A$ is unbounded on the right (resp.~left).
In \hyperref[sec:allOAGareAI]{Section 4}, we prove our 
main result: any non-trivial ordered abelian group is augmentable by infinites (Theorem 
\ref{thm:strongleftaugmentabilityOAG}), and discuss how to tell that the quotient formed by the 
infinites is ($p$-)divisible (Corollary \ref{cor_div}). As an application of our main result, we 
discuss consequences for henselian valued fields in \hyperref[sec:def]{Section 5}, in particular
Proposition \ref{power} and Theorem \ref{def_fieldwise}. 
In the \hyperref[sec:app]{appendix},
we give an alternative proof of characterizing which coloured multi-orders are augmentable (Corollary \ref{cor_iksrat}).

\subsection*{Acknowledgements}
We deeply thank Salma Kuhlmann for pointing out to us the work of Delon and Lucas \cite{DelonLucas}. We thank Sylvy Anscombe, Philip Dittmann, Arno Fehm, Immanuel Halup\-czok, Mathias Stout, and Floris Vermeulen for their useful comments and discussions.

\section{Spines}\label{sec:spines}

We recall the language $\mathcal{L}_{\text{syn}}$ introduced by Cluckers and Halupczok in \cite{CH11}. We begin by describing the \emph{auxiliary sorts} $\mathcal{S}_n$, $\mathcal{T}_n$ and $\mathcal{T}_n^+$ for each $n \in \mathbb{N}$, $n>0$.
\begin{definition} \label{DefinitionAuxiliarySorts}
Fix a natural number $n>0$.
\begin{enumerate}
    \item For $g \in G \setminus nG$, let $G_g^n$ be the largest convex subgroup $H$ of $G$ such that $g \notin H+nG$.
    
    Define $\mathcal{S}_n := (G\setminus nG) / {\sim}$, with $g \sim g'$ if and only if $G_g^n=G_{g'}^n$, and let $\mathfrak{s}_n \colon G \twoheadrightarrow \mathcal{S}_n$ be the canonical map. Denote by $G_{\alpha}$ the convex subgroup $G_g^n$, with $\alpha=\mathfrak{s}_n(g)$.
    
    \item For $g \in G$, set $H_g^n=\bigcup_{h \in G, g \notin G_h^n}G_h^n$ with the convention that a union over the empty set is $\{0\}$.
   
    Define $\mathcal{T}_n := G / {\sim}$, with $g \sim g'$ if and only if $H_g^n=H_{g'}^n$, and let $\mathfrak{t}_n \colon G \twoheadrightarrow \mathcal{T}_n$ be the canonical map. Denote by $G_{\alpha}$ the convex subgroup $H_g^n$, with $\alpha=\mathfrak{t}_n(g)$.
    
    \item For $g\in G$, let $H_{g^+}^n=\bigcap_{h \in G, g \in G_h^n}G_h^n$ , where the intersection over the empty set is $G$.
    
    Define $\mathcal{T}_n^+$ the quotient $G/{\sim}$, with $g \sim g'$ if and only if $H_{g^+}^n=H_{{g'}^+}^n$. Let $\mathfrak{t}_n^+$ the natural map and for $g\in G$, $\mathfrak{t}_n^+(g)=:\beta^+\in \Tcal_n^+$, denote $H_{g^+}$ by $G_{\beta^+}^n$. Note that $ \mathcal{T}_n^+$ is a copy of $\mathcal{T}_n$, with potentially an additional initial point:  the map $\iota^+:\Tcal_n \rightarrow \Tcal_n^+, \mathfrak{t}_n(g) \mapsto \mathfrak{t}_n^+(g)$ is injective and decreasing, with regard to the ordering defined in \Cref{DefinitionLsyn}. Now $\Tcal_n^+ \setminus \iota^+( \Tcal_n)$ is either empty or the first point of $\Tcal_n^+$.
\end{enumerate}

\end{definition}
   
In \cite{CH11}, it is proved that, for each $n>0$, the convex subgroups of the three families in the above definition $\{G_{\alpha}\}_{\alpha \in \mathcal{S}_n}, \{G_{\alpha}\}_{\alpha \in \mathcal{T}_n}$ and $ \{G_{\alpha}\}_{\alpha \in \mathcal{T}_n^+}$ are uniformly definable in  $\mathcal{L}_{\text{oag}}=\{0,+,-,<\}$. It follows that in any theory of ordered abelian groups, all the auxiliary sorts are imaginary sorts of $\mathcal{L}_{\text{oag}}$.

For any $\alpha \in \bigcup_{n \in \mathbb{N},n>0} (\Scal_n\cup\Tcal_n\cup\Tcal_n^+)$ and $m \in \mathbb{N}, m>0$, we also set
\begin{equation}
G_{\alpha}^{[m]} := \bigcap_{\phantom{MM}\mathclap{\substack{G_\alpha\subsetneq H\subseteq G\\H\text{ convex subgroup}}}\phantom{MM}} (H + mG).
\end{equation}
Now we are able to present the complete definition of $\mathcal{L}_{\text{syn}}$. 
\begin{definition}\label{DefinitionLsyn}
The language $\mathcal{L}_{\text{syn}}$ consists of the following:
\begin{enumerate}[label=(\alph*)]
    \item\label{a} The main sort $(G,0,+,-,<,(\equiv_m)_{m \in \mathbb{N}})$;
    \item\label{b} the auxiliary sorts $\mathcal{S}_p,\mathcal{T}_p, \mathcal{T}_p^+$, for each prime $p$, with the binary relations $\leqslant_{p,q}$ on ($\mathcal{S}_p \cup \mathcal{T}_p \cup \mathcal{T}_p^+) \times (\mathcal{S}_q \cup \mathcal{T}_q \cup \mathcal{T}_q^+)$ (each pair of primes $(p,q)$ giving rise to nine binary relations), often simply denoted by $\leqslant$, defined by $\alpha \leqslant \alpha'$ if and only if $G_{\alpha'} \subseteq G_{\alpha}$\footnote{In \cite{CH11}, the order is defined so that it corresponds to the inclusion of convex subgroups. Here, we use the reverse order, so that $\mathfrak{t}_p$ is a pre-valuation.};
    \item\label{c} the canonical maps $\mathfrak{s}_p \colon G \twoheadrightarrow \mathcal{S}_p$, $\mathfrak{t}_p \colon G \twoheadrightarrow \mathcal{T}_p$ and $\mathfrak{t}_p^+ \colon G \twoheadrightarrow \mathcal{T}_p^+$, for each prime $p$;

    \item\label{d} a unary predicate $x=_\bullet k_\bullet$ on $G$, for each $k \in \Zbb \setminus \{0\}$, defined by $g=_\bullet k_\bullet$ if and only if there exists a convex subgroup $H$ of $G$ such that $G/H$ is discrete and $g \mod H$ is equal to $k$ times the smallest positive element of $G/H$, for every $g \in G$; 
    \item\label{e} a unary predicate $x\equiv_{\bullet m} k_\bullet$ on $G$, for each $m \in \mathbb{N} \setminus \{0\}$ and $k \in \{1,\dots,m-1\}$, defined by $g \equiv_{\bullet m} k_\bullet$ if and only if there exists a convex subgroup $H$ of $G$ such that $G/H$ is discrete and $g \mod H$ is congruent modulo $m$ to $k$ times the smallest positive element of $G/H$, for every $g \in G$;
    \item\label{f} a unary predicate $D_{p^r}^{[p^s]}(x)$ on $G$, for each prime $p$ and each $r,s \in \mathbb{N} \setminus \{0\}$ with $s \geqslant r$, defined by $D_{p^r}^{[p^s]}(g)$ if and only if there exists an $\alpha \in \mathcal{S}_p$ such that $g \in G_{\alpha}^{[p^s]}+p^rG$ and $g \notin G_{\alpha}+p^rG$, for every $g \in G$;
    \item \label{g}a unary predicate discr$(x)$ on the sort $\mathcal{S}_p$, with $p$ prime, defined by discr$(\alpha)$ if and only if $G/G_{\alpha}$ is discrete, for every $\alpha \in \mathcal{S}_p$;
    \item\label{h} unary $P^{[p^n]}_{l},Q^{p}_l$ predicates on the sort $\mathcal{S}_p$, with $p$ prime, for each $l,n \in \mathbb{N} \setminus \{0\}$, defining the sets
    \begin{gather*}
    P^{[p^n]}_{l}:=\{\alpha \in \mathcal{S}_p \mid \dim_{\mathbb{F}_p}(G_{\alpha}^{[p^n]}+pG)/(G_{\alpha}^{[p^{n+1}]}+pG)=l\} \text{ and} \\
    Q^{p}_l:= \{\alpha \in \mathcal{S}_p \mid \dim_{\mathbb{F}_p}(G_{\alpha}^{[p]}+pG)/(G_{\alpha}+pG)=l\}\footnotemark.\end{gather*}
\footnotetext{ In \cite{CH11}, for any $N$, a predicate interpreting $\{ \alpha \mid \dim_{\mathbb{F}_p}(G_{\alpha}^{[p^N]}+pG)/(G_{\alpha}+pG)=l\}$ is added to the language. It can be recovered from the predicates above since $\dim_{\mathbb{F}_p}(G_{\alpha}^{[p^N]}+pG)/(G_{\alpha}+pG)=\dim_{\mathbb{F}_p}(G_{\alpha}^{[p]}+pG)/(G_{\alpha}+pG) - \sum_{n<N} \dim_{\mathbb{F}_p}(G_{\alpha}^{[p^n]}+pG)/(G_{\alpha}^{[p^{n+1}]}+pG) $. }
\end{enumerate}  
We denote by $\mathcal{A}_{p}$ the set of sorts $\Scal_p,\Tcal_p, \Tcal_p^+$ for any fixed $p$, equipped with the orders and the unary predicates discr, $P^{[p^n]}_{l}$ and $Q^p_{l}$. We call the set $\Acal_{p}$ the $p$-\emph{spine} of $G$. The \emph{spines} of $G$, denoted by $\Acal$, are then defined as the set of sorts $\Acal_{p}$ for all primes $p$ equipped furthermore with the orderings $\leqslant_{p,q}$ for all pairs of primes $(p,q)$.
\end{definition}

\begin{fact}[{\cite[Theorem 1.13]{CH11}}]\label{FactQUantifierEliminationLsyn}
The $\mathcal{L}_{syn}$-theory of ordered abelian groups eliminates quantifiers relative to the spines.
\end{fact}

    Since $\Acal$ is a closed set of sorts, it is stably embedded and pure in the sense of \cite[Appendix A]{Rid17}.
    This result of quantifier elimination says that any $\mathcal{L}_{\text{syn}}$-formula $\phi(\bar{x},\bar{\eta})$, with $G$-variables $\bar{x}$ and $\mathcal{A}$-variables $\bar{\eta}$, is a boolean combination of formulas of the form
\begin{itemize}
	\item $\psi(\bar{x})$, where $\psi$ is quantifier-free in the language 
 \[(G,0,+,-,<,(\equiv_m)_{m \in \mathbb{N}}, (=_\bullet k_\bullet)_{k\in \Nbb}, (\equiv_{\bullet m} k_\bullet)_{m,k\in \Nbb}, (D^{[p^s]}_{p^r})_{p\in \mathbb{P}, s\geqslant r \in \Nbb }),\]
  and
	\item $\chi(\bar{x},\bar{\eta}):=\xi((\mathfrak{s}_{p_j}(\sum_{i<n}z_i^jx_i),\mathfrak{t}_{p_j}(\sum_{i<n}z_i^jx_i))_{z_0^j, \dots, z_{n-1}^j \in \Zbb, j\in J},\bar{\eta})$, where $\xi$ is an $\mathcal{L}_{\text{syn}\vert_{\mathcal{A}}}$-formula, and $J$ is a set of primes (with potential repetitions).
\end{itemize}


\newcommand{\Lsyn}{\Lcal_{\mathrm{syn}}}
\newcommand{\Loag}{\Lcal_{\mathrm{oag}}}

Note that since (the interpretation of) every new symbol in $\Lsyn$ is definable, for any two ordered abelian groups, we have $G\equiv G'$ in $\Loag$ if and only if $G\equiv G'$ in $\Lsyn$, and alike for $G\substruct G'$. In particular, if $G\substruct G'$, then $\Acal(G)\substruct \Acal(G')$. Furthermore, by quantifier elimination, we have:

\begin{fact}[CH/S transfer principle for ordered abelian groups]\label{fact:Schmitt}\label{fact:strongSchmitt}
Let $G$ and $G'$ be two ordered abelian groups. Then:
\[G\equiv G' \ \text{ if and only if } \Acal(G) \equiv \Acal(G')\ .\]
Furthermore, if $G$ is an $\Lsyn$-substructure of $G'$, then:
$$G\substruct G' \ \text{ if and only if } \Acal(G) \substruct \Acal(G')\ .$$
\end{fact}

We call this fact the ``CH/S transfer principle'' as an echo to the AK/E transfer principle, after Schmitt, who stated in very different terms a similar transfer principle in his habilitation thesis \cite{Sch82}, and Cluckers-Halupczok, who derived the language $\Lsyn$ for quantifier elimination \cite{CH11}, giving us this more refined version of the transfer principle.

\ 

In general, given $G'$ an ordered abelian group and $G\leq G'$ a subgroup, 
we do not have in general an inclusion of spines, that is, it might be that $\Acal(G)$ is not a subset of $\Acal(G')$. For example, $\mathbb Z\subseteq\mathbb Q$, but $\Scal_p(\mathbb Z)$ has 1 point (for $\{0\}$) where $\Scal_p(\mathbb Q)$ is empty.

We will now derive some conditions under which (convex) subgroups can be seen as $\Lsyn$-substructures.
\newcommand{\sfp}{\mathfrak s_p}
\newcommand{\tfp}{\mathfrak t_p}
\newcommand{\tfpp}{\mathfrak t_p^+}

\begin{lemma}\label{lem:witness-of-cvx}
 Let $G$ be an ordered abelian group and $H\subseteq G$ be a convex subgroup of $G$. Consider $G$ as an 
 $\Lsyn$-structure, specifically with the maps $(\sfp,\tfp,\tfpp)$. Then convex subgroups of $H$ which are $G$-definable in $G$ are $H$-definable in $G$, in the following sense:
\begin{enumerate}
 \item\label{s} Let $a\in G$ such that $G_{\mathfrak s_p(a)}\subsetneq H$. Then there exists $b\in H$ such that $\mathfrak s_p(b)=\mathfrak s_p(a)$.
 \item\label{t} Let $a\in G$ such that $G_{\mathfrak t_p(a)}\subsetneq H$. Then there exists $b\in H$ such that $\mathfrak t_p(b)=\mathfrak t_p(a)$.
 \item\label{t+} Let $a\in G$ such that $G_{\mathfrak t_p^+(a)}\subseteq H$. Then $a\in H$.
 \end{enumerate}
\end{lemma}
\begin{proof}
We first prove (\ref{s}). Recall that $G_{\mathfrak s_p(a)}$ is the largest subgroup of $G$ such that $a$ is not $p$-divisible modulo $G_{\mathfrak s_p(a)}$. Since $H$ is strictly larger than $G_{\mathfrak s_p(a)}$, $a$ is $p$-divisible modulo $H$, that is: $a=pa'+b$ for some $a'\in G$ and $b\in H$. Now $b=a-pa'$ is $p$-divisible modulo some subgroup $\Delta$ iff $a$ is $p$-divisible modulo $\Delta$, that is, $\mathfrak s_p(a)=\mathfrak s_p(b)$.

For (\ref{t}), since $G_{\mathfrak t_p(a)}\subsetneq H$, take $b\in H\setminus G_{\mathfrak t_p(a)}$. By definition, $G_{\mathfrak t_p(a)}\subseteq G_{\mathfrak t_p(b)}$. If this inclusion is strict, then $a\in G_{\mathfrak t_p(b)}\subseteq H$. If not, then $\mathfrak t_p(a)=\mathfrak t_p(b)$ as wanted.

For (\ref{t+}), it is immediate, as by definition $a\in G_{\mathfrak t_p^+(a)}\subseteq H$.
\end{proof}

\begin{lemma}\label{lem:cvx-are-ok}
 Let $G$ be an ordered abelian group and $H\subseteq G$ be a convex subgroup of $G$. We consider $H$ as an ordered abelian group on its own, and we equip $G$ and $H$ with $\Lsyn$, specifically, with maps $(\sfp,\tfp,\tfpp)\colon H\rightarrow \Acal(H)$ and $(\sfp',\tfp',{\tfpp}')\colon G\rightarrow \Acal(G)$. Then these maps agree on proper convex subgroups of $H$, that is:
 
 \begin{enumerate}
  \item let $a\in H\setminus pH$ and let $\alpha=\mathfrak s_p(a)\in\Acal(H)$ and $\alpha'=\mathfrak s_p'(a)\in\Acal(G)$. Then $G_\alpha$, as a subgroup of $H$, is equal to $G_{\alpha'}$, as a subgroup of $G$.
  \item let $a\in H$ and let $\alpha=\mathfrak t_p(a)\in\Acal(H)$ and $\alpha'=\mathfrak t_p'(a)\in\Acal(G)$. Then $G_\alpha$, as a subgroup of $H$, is equal to $G_{\alpha'}$, as a subgroup of $G$. 
  \item let $a\in H$ and let $\alpha=\mathfrak t^+_p(a)\in\Acal(H)$ and $\alpha'=\mathfrak {t^+_p}'(a)\in\Acal(G)$.
  \begin{enumerate}
  \item If $G_{\alpha}\neq H$, then $G_\alpha$, as a subgroup of $H$, is equal to $G_{\alpha'}$, as a subgroup of $G$.
  \item If $G_{\alpha}=H$, then $G_{\alpha'}\supseteq H$.
  \end{enumerate}
 \end{enumerate}
\end{lemma}
\begin{proof}
 \begin{enumerate}
  \item Let $a\in H\setminus pH$. Then also $a\in G\setminus pG$ since $H$ is a convex subgroup of $G$. Let $\alpha=\mathfrak s_p(a)\in\Acal(H)$ and $\alpha'=\mathfrak s_p'(a)\in\Acal(G)$. Then by definition, $G_\alpha$ is the largest convex subgroup of $H$ such that $a$ is not $p$-divisible modulo $G_\alpha$. In particular, $G_\alpha$ is a subgroup of $H$ and thus also of $G$. Since $G_{\alpha'}$ is the largest convex subgroup of $G$ such that $a$ is not $p$-divisible in modulo $G_{\alpha'}$, in particular, $G_{\alpha'}$ is larger than $G_{\alpha}$. But since $G_{\alpha'}\subseteq[-|a|,|a|]\subseteq H$, $G_{\alpha'}$ is itself a convex subgroup of $H$ such that $a$ is not $p$-divisible modulo $G_{\alpha'}$, so $G_{\alpha}$ is larger than $G_{\alpha'}$, that is, they are equal.
  
  \item Let $a\in H$ and let $\alpha=\mathfrak t_p(a)\in\Acal(H)$ and $\alpha'=\mathfrak t_p'(a)\in\Acal(G)$. Then by definition:
  $$G_\alpha=\bigcup_{\phantom{M}\mathclap{\substack{b\in H\\a\notin G_{\mathfrak s_p(b)}}}\phantom{M}}G_{\mathfrak s_p(b)}\;\text{ and }\;G_{\alpha'}=\bigcup_{\phantom{M}\mathclap{\substack{b\in G\\a\notin G_{\mathfrak s_p'(b)}}}\phantom{M}}G_{\mathfrak s_p'(b)}.$$
  Clearly $G_{\alpha}\subseteq G_{\alpha'}$, since the union is taken over a larger set of indices, and since by the proof of (1) above $\sfp$ and $\sfp'$ agree on $H$. Now assume, towards a contradiction, that this inclusion is strict. This means that there is $b\in G\setminus pG$ such that $a\notin G_{\mathfrak s_p'(b)}$ and $G_{\mathfrak t_p(a)}\subsetneq G_{\mathfrak s_p'(b)}$. By \Cref{lem:witness-of-cvx}, we can assume $b\in H$. But then $G_{\mathfrak s_p'(b)}=G_{\sfp(b)}$ appears in the definition of $G_{\mathfrak t_p(a)}$, that is, $G_{\sfp(b)}\subseteq G_{\mathfrak t_p(a)}$, a contradiction.
  
  \item Let $a\in H$ and let $\alpha=\mathfrak t^+_p(a)\in\Acal(H)$ and $\alpha'=\mathfrak {t^+_p}'(a)\in\Acal(G)$. By definition, $G_\alpha=\bigcap_{b\in H, a\in G_{\mathfrak s_p(b)}}G_{\mathfrak s_p(b)}$ and $G_{\alpha'}=\bigcap_{b\in G, a\in G_{\mathfrak s_p'(b)}}G_{\mathfrak s_p'(b)}$.
  \begin{enumerate}
  \item If $G_{\alpha}\neq H$, then the set of indices is not empty. Since $\sfp$ and $\sfp'$ agree on $H$, and since any subgroup of $H$ can be defined in $H$ by \Cref{lem:witness-of-cvx}, we have:
  $$G_{\alpha'}=\bigcap_{\phantom{M}\mathclap{\substack{b\in G\\a\in G_{\mathfrak s_p'(b)}}}\phantom{M}}G_{\mathfrak s_p'(b)}=\bigcap_{\phantom{M}\mathclap{\substack{b\in H\\ a\in G_{\mathfrak s_p(b)}}}\phantom{M}}G_{\mathfrak s_p(b)}=G_{\alpha}.$$
  
  \item If $G_{\alpha}=H$, then the intersection defining $G_{\alpha}$ was over the emptyset, that is, there is no convex subgroup of $H$ definable in $H$ and containing $a$. But then there is also no convex subgroup of $H$ definable in $G$ containing $a$, and so any definable convex subgroup containing $a$ must be larger than $H$, that is, $H\subseteq G_{\alpha'}$.
  \end{enumerate}
 \end{enumerate}
\end{proof}

\newcommand{\Sp}{\Scal_p}
\newcommand{\Tp}{\Tcal_p}
\newcommand{\Tpp}{\Tcal_p^+}
\begin{corollary}\label{cor:cvx-are-substruct}
 Let $G$ be an ordered abelian group and $H$ be a convex subgroup of $G$. We equip $G$ and $H$ with the language $\Lsyn$ independently of one another. 
 Assume the following holds:
 \begin{equation} \text{For all primes $p$ and for all $a\in H$, if }G_{{\mathfrak{t}_p^+}^H(a)}= H \text{ then } G_{{\mathfrak{t}_p^+}^G(a)}= G. \tag{$\star$} \label{condition:CFirst} \end{equation}
 Equivalently, for all prime $p$ such that there is $a\in H$ with $G_{{\mathfrak{t}_p^+}^H(a)}= H$, $G/H$ is divisible by $p$. 
 Then $H$ is naturally an $\Lsyn$-substructure of $G$. More precisely, we define an embedding ${\iota\colon H\rightarrow G}$ extending the inclusion $H\subseteq G$ by setting, for any $a\in H$ and any prime $p$, $\iota(\sfp^H(a))=\sfp^G(a)$, $\iota(\tfp^H(a))=\tfp^G(a)$ and $\iota({\tfpp}^H(a))={\tfpp}^G(a)$. In this way, $\iota(H)$ is the $\Lsyn$-structure obtained by the restricting the symbols of $\Lsyn$ from $G$ to $H$.
 \end{corollary}

 \begin{example}
     The condition (\ref{condition:CFirst}) cannot be dropped. For example, consider the groups $G\coloneqq \mathbb{Z}_{(2)} \oplus \mathbb{Z}$ and $H \coloneqq \{0\} \oplus \mathbb{Z}$. The spines of $H$ do not embed in the spines of $G$ precisely because  $G_{{\mathfrak{t}_2^+}^G((0,1))}= H$ and $G_{{\mathfrak{t}_3^+}^G((0,1))}= G$, while $G_{{\mathfrak{t}_2^+}^H((0,1))}=G_{{\mathfrak{t}_3^+}^H((0,1))}=H$.
 \end{example}

\begin{proof}[Proof of \Cref{cor:cvx-are-substruct}]
 The fact that $\iota$ is well-defined follows from \Cref{lem:cvx-are-ok} and \Cref{lem:witness-of-cvx}. What remains to check that $\iota$ is an $\Lsyn$-embedding: 
  All computations of proper convex subgroups of $H$ yield the same convex subgroup when computed in $G$ (cf. \Cref{lem:cvx-are-ok}).
 The only potential obstruction comes when a point in $\mathcal{T}_p^+$ defines $H$ (which will then necessarily be an initial point in $\Tpp$). Then condition (\ref{condition:CFirst})
 ensures that doing the same computation in $G$ results with $G$
(rather than a proper convex subgroup containing $H$).
\end{proof}

We can even be more precise: since points in the spines of $H$ do not change when computed in $G$ -- with the possible exception of initial points in $\Tpp$ -- we have the following:

\begin{corollary}[See {\cite[Lemma 1-1]{DelonLucas}}]\label{cor:spines-of-cvx}
Let $G$ be an ordered abelian group and $H$ a convex subgroup of $G$. We consider $G$ and $H$ each
as $\Lsyn$-structures. Let $\Ccal=\sset{\alpha\in\Acal(H)}{G_\alpha=H}$. Then $\Acal(H)\setminus\Ccal$ is an end segment of $\Acal(G)$.
\end{corollary}

\begin{proof}
    Note that $H$ embeds into $G$ in the restricted language $\Lsyn^{-}$ where the sorts
    $\Tpp$ (and corresponding functions) are dropped. 
    Any point in $\Tpp$ which is not initial 
    is already in $\Tp$. The set $\Ccal$
    contains all initial points of $\Tpp$, and so $\Acal(H)\setminus\Ccal$ naturally embeds into $\Acal(G)$. 
\end{proof}

\section{Coloured multi-orders}\label{sec:multiord}

In the section, we study spines as abstract structures, detached from the ordered abelian group they come from. As such, they are coloured multi-orders (see \Cref{def:multiord}), that is, a many-sorted structure where each sort is a linear order, equipped furthermore with pre-orders between distinct sorts, as well as unary predicates. 

We prove notably that given a multi-order $A$ and an elementary extension $A^*$ of $A$, the convex hull of $A$ in $A^*$, the left hull and the right hull of $A$ in $A^*$ are intermediate elementary extensions, and that a multi-order is augmentable if and only if one of its components is unbounded.

\begin{definition}\label{def:multiord}
Let $P$ be a set. The language of $P$-multi-orders consists of one sort $A_p$ for each $p\in P$ and for each pair $(p,q)\in P^2$ (including pairs of the form $(p,p)$) a binary relation $\leqslant_{p,q}$ on $A_p\times A_q$. A $P$-\emph{multi-order} is a structure in this language such that for each $p\in P$, $\leqslant_{p,p}$ is a linear order on each component $A_p$, and $\bigcup_{p,q\in P}\leqslant_{p,q}$ is a linear pre-order on $\bigcup_{p\in P}A_p$.

A coloured $P$-multi-order is a $P$-multi-order augmented by arbitrarily many unary predicates.
\end{definition}

By linear preorder, we mean that two elements $a$ and $b$, respectively in $A_p$ and $A_q$ for $p\neq q$, are always comparable by $\leqslant_{p,q}$ and by $\leqslant_{q,p}$; however, it might happen that both $a\leqslant_{p,q}b$ and $b\leqslant_{q,p}a$ hold, but still $a\neq b$ since they do not lie in the same sort.

Note that spines of ordered abelian groups can naturally
be seen as coloured $3\mathbb P$-multi-orders where $3\mathbb P$ denotes three disjoint copies of the set of primes, equipped with orderings $\leqslant_{p,q}$ and colours $\discr$, $P^{[p^n]}_l$ and $Q^p_l$ as defined in \Cref{DefinitionAuxiliarySorts}.

When needed, we will often write $x\leqslant y$ without specifying for which $p,q$ the relation $\leqslant_{p,q}$ is considered.

\begin{definition}
 Let $(I,<)$ be a linear order and for each $i\in I$ let $A^i$ be a coloured $P$-multi-order (each of them with the same set of colours). The sum $A=\sum_{i\in I}A^i$ is a coloured $P$-multi-order defined by:
 \begin{itemize}
 \item for each $p\in P$, $A_p=\bigsqcup_{i\in I}A^i_p$;
 \item for each $p,q\in P$, for each $a\in A_p$ and each $b\in A_q$, $A\satisf a\leqslant_{p,q}b$ iff either both $a$ and $b$ lie in the same $A^i$ and $A^i\satisf a\leqslant_{p,q} b$, or $a\in A^i_p$, $b\in A^j_q$ and $I\satisf i<j$.
 \item for each unary predicate $C$ and each $a\in A$, $A\satisf C(a)$ iff $A^i\satisf C(a)$, where $a\in A^i$. 
 \end{itemize}
\end{definition}

Intuitively, $A+B$ corresponds to a copy of $A$ followed by a copy of $B$. Instead of writing $\sum_{i\in I}A^i$, we will often write $A=A^1+\cdots+A^n$.

The following lemma is an immediate generalization to multi-orders of \cite[Theorem 2.1]{Rub74}, also obtained using Ehrenfeucht-Fraïssé games:
\begin{lemma}\label{lem:sumofmo}
Let $(I,<)$ be a linear order and for each $i\in I$ let $A^i$, $B^i$ be coloured $P$-multi-orders (all with the same set of colours), and write $A=\sum_{i\in I} A^i$ and $B=\sum_{i\in I}B^i$.
\begin{enumerate}
\item If $A^i\equiv B^i$ for all $i\in I$, then $A\equiv B$.
\item If $A^i\substruct B^i$ for all $i\in I$, then $A\substruct B$.
\end{enumerate}
\end{lemma}

\begin{definition}
 Let $B$ be a coloured $P$-multi-order and let $A\subseteq B$. We define:
 \begin{enumerate}
  \item $\langle A\rangle_B=\sset{x\in B}{\exists a,a'\in A(a\leqslant x\leqslant a')}$, the convex hull of $A$ in $B$,
  \item $B_{<A}=\sset{x\in B}{\forall a\in A(x<a)}$,
  \item $B_{>A}=\sset{x\in B}{\forall a\in A(x>a)}$,
  \item $B_{\leqslant A}=B_{<A}+\langle A\rangle_B$, the left hull of $A$ in $B$,
  \item $B_{\geqslant A}=\langle A\rangle_B+B_{>A}$, the right hull of $A$ in $B$.
 \end{enumerate}
\end{definition}

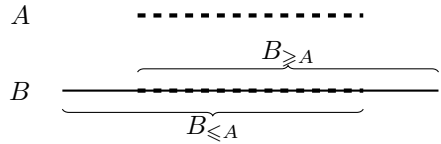
    \begin{center}
    \begin{tikzpicture}
    \draw (-0.3,1) node[anchor=east] {$B$}; 
    \draw[dashed,ultra thick] (1,2) -- (4,2);
    \draw[decoration={brace,raise=7pt},decorate]
  (1,1) -- (5,1) ;
    \draw (3,1.2) node[anchor=south] {$B_{\geqslant A}$};
    \draw[decoration={brace,raise=7pt},decorate]
  (4,1) -- (0,1) ;
    \draw (2,0.2) node[anchor=south] {$B_{\leqslant A}$};
    \draw (-0.3,2) node[anchor=east] {$A$};
    \draw[thick] (0,1) -- (5,1);
    \draw[dashed,ultra thick] (1,1) -- (4,1);
\end{tikzpicture}
\end{center}

\begin{lemma}[See {\cite[Example 1]{DelonLucas}}]\label{lem:hulls!}
 Let $A,B$ be coloured $P$-multi-orders with $A\substruct B$. Then:
 \begin{enumerate}
  \item $A\substruct \langle A\rangle_B\substruct B$,
  \item $A\substruct B_{\leqslant A}\substruct B$, and
  \item $A\substruct B_{\geqslant A}\substruct B$.
 \end{enumerate}
\end{lemma}

\begin{proof}
 Let $a\in A$. With the notation above, we have $B_{\leqslant a}=\sset{x\in B}{x\leqslant a}$ and $A_{>a}=\sset{x\in A}{x>a}$. Denote $B_{\leqslant a}+A_{>a}$ by $C_a$. We have naturally $A_{\leqslant a}\substruct B_{\leqslant a}$, therefore by \Cref{lem:sumofmo}, $A=A_{\leqslant a}+A_{>a}\substruct C_a\substruct B$. The chain of structures $(C_a)_{a\in A}$ is thus elementary, and by Tarski's chain lemma \cite[Theorem~2.1.4 \& Exercise~2.1.1]{TZ12}, $A\substruct \bigcup_{a\in A} C_a\substruct B$, and $\bigcup_{a\in A} C_a$ is exactly $B_{\leqslant A}$.
 
 Reversing the ordering gives a proof of $A\substruct B_{\geqslant A}\substruct B$. For the convex hull, we simply note that ${(B_{\leqslant A})}_{\geqslant A}=\langle A\rangle_B$.
\end{proof}

We now use this fact to characterize ``augmentable'' multi-orders:

\begin{definition}
 Let $X$ and $Y$ be coloured $P$-multi-orders. We say that $X$ is augmentable on the right (resp.~on the left) by $Y$ if $X\substruct X+Y$ (resp.~$X\substruct Y+X$).
 
 We say that a coloured $P$-multi-order augmentable on the right (resp.~on the left) if there exists a non-empty $P$-multi-order $Y$ such that $X$ is augmentable on the right (resp.~on the left) by $Y$.
\end{definition}

\begin{definition}
 Let $A$ be a coloured $P$-multi-order and let $p\in P$. We say that a component $A_p$ is bounded on the right (resp.~on the left) if there exists $a\in A$, say $a\in A_q$, such that $x\leqslant_{p,q}a$ (resp.~$a\leqslant_{q,p}x$) for all $x\in A_p$. Otherwise, we say that $A_p$ is unbounded on the right (resp.~on the left).
\end{definition}

\begin{lemma}\label{lem:aug-mo}
 Let $A$ be a coloured $P$-multi-order. Then $A$ is augmentable on the right (resp.~on the left) iff for some $p\in P$, $A_p$ is unbounded on the right (resp.~on the left).
\end{lemma}

\begin{proof}
We do the proof on the right, reversing the ordering gives a proof on the left. First, suppose that every $A_p$ is bounded on the right, that is, for every $p\in P$, there is $a\in A_q$ for some $q\in P$ such that $\forall x\in A_p$, $x\leqslant_{p,q} a$. Let $X$ be a coloured $P$-multi-order such that $A\substruct A+X$. We aim to prove that $X$ is empty. Indeed, we have $A_p\substruct (A+X)_p\cong A_p+X_p$; but ``$A_p$ is bounded by $a$'' is a first-order formula with parameter $a\in A$, thus it is also satisfied by $A+X$ and $A_p+X_p$ must also be bounded by $a$. But by definition of $A+X$ any element of $X_p$ is strictly bigger than any element of $A$, in particular bigger than $a$; hence $X_p$ must be empty. Since this holds for all $p\in P$, $X$ must be empty, and we conclude that $A$ is not augmentable on the right.

Conversely, suppose that some $A_p$ is unbounded. This means that for any $a\in A$, say $a\in A_q$, there exists $b\in A_p$ such that $a\leqslant_{q,p}b$. The following ``type at infinity'' is thus satisfiable:
$$\pi_{\infty}(x)=\sset{a\leqslant_{q,p} x}{a\in A_q}_{q\in P}$$

Let $B\superstruct A$ be $|A|^+$-saturated. Then in particular, $B$ contains a realisation of $\pi_{\infty}$, that is, $B_{>A}\neq\emptyset$.

Now by \Cref{lem:hulls!}, we have $A\substruct B_{\leqslant A}$. By \Cref{lem:sumofmo}, we have $A+X\substruct B_{\leqslant A}+X$ for any coloured $P$-multi-order $X$; in particular, $A+B_{>A}\substruct B_{\leqslant A}+B_{>A}=B$. Now $A\substruct B$ and $A\subset A+B_{>A}\substruct B$, therefore, $A\substruct A+B_{>A}$.
\end{proof}

\begin{corollary}
 Let $A$ be a coloured $P$-multi-order. Let $C$ be the set of cuts of $A$ which are satisfiable but not realized in $A$. Let $C'\subset C$, then there is an elementary extension of $A$ where every cut in $C'$ is realized and every cut in $C\setminus C'$ is not.
\end{corollary}

\begin{proof}
 Given a cut $(A_{<x}$,$A_{>x})\in C$, since it is satisfiable, either $A_{<x}$ is unbounded on the right or $A_{>x}$ is unbounded on the left. In any case, applying \Cref{lem:aug-mo}, there is $X$ such that $A\substruct A_{<x}+X+A_{>x}$; that is, there is an elementary extension of $A$ realizing only this cut.
 
 Now fix $C'\subseteq C$ and let $\lambda=|C'|$ enumerate $C'$. Fix $A_0=A$ and for each $i<\lambda$, let $A_{i+1}$ be an elementary extension of $A_i$ realizing only the cut number $i$ as explained above. For limit ordinals $i<\lambda$, let $A_i=\bigcup_{j<i}A_j$. Now by Tarski's chain lemma \cite[Theorem 2.1.4]{TZ12} $A_\lambda=\bigcup_{i<\lambda} A_i$ is an elementary extension of $A$, and it realizes exactly the cuts in $C'$.
\end{proof}

\section{Augmentability by infinites}\label{sec:allOAGareAI}
\begin{lemma}[See {\cite[Proposition 1]{DelonLucas}}]\label{lem:ConvxHull}
    Let $G$ be an ordered abelian group and $G'$ an elementary extension. Then the convex hull $H:=\langle G \rangle$ of $G$ in $G'$ is also an elementary extension of $G$ and an elementary substructure of $G'$.
\end{lemma}
\begin{proof}
Since $G\substruct G'$, clearly $\Acal(G)\substruct \Acal(G')$. 
\begin{claim} 
The pair $G'$ and $H$ satisfies (\ref{condition:CFirst}) from \Cref{cor:cvx-are-substruct}. 
\end{claim}
\begin{proof}[Proof of Claim 1]
Assume $G_{{\mathfrak t_p^+}^H(a)}=H$ for some prime $p$ and some $a\in H$. 
Since $H$ is the convex hull of $G$, there exists $b\in G$ such that 
$|b|>|a|$.

But now $G_{{\mathfrak t_p^+}^G(b)}=G$: indeed, $G_{{\mathfrak t_p^+}^H(b)}=H$, then by \Cref{lem:cvx-are-ok}, $G_{{\mathfrak t_p^+}^{G'}(b)}\supseteq H$, and since $G\subseteq H$, we have $G_{{\mathfrak t_p^+}^{G'}(b)}\supseteq G$. By elementarity, this can only happen if $G_{{\mathfrak t_p^+}^{G}(b)}=G$.

Hence $G_{{\mathfrak t_p^+}^{G'}(b)}=G'$, using elementarity once again.

Now $a\in G_{{\mathfrak t_p^+}^{G'}(b)}$ and $b\in H=G_{{\mathfrak t_p^+}^{H}(a)}\subseteq G_{{\mathfrak t_p^+}^{G'}(a)}$, thus we must have ${\mathfrak t_p^+}^{G'}(b)={\mathfrak t_p^+}^{G'}(a)$ and finally $G_{{\mathfrak t_p^+}^{G'}(a)}=G'$, proving that the pair $H$ and $G'$ satisfies condition (\ref{condition:CFirst}). 
\end{proof}

Therefore, by \Cref{cor:cvx-are-substruct}, $H$ is an $\Lsyn$-substructure of $G'$.
    
    As in \Cref{cor:spines-of-cvx}, let $$\Ccal=\sset{\alpha\in\Acal(G)}{G_\alpha=G}.$$ Since ``$G_\alpha=G$'' is a first order property, these points are also in $\Acal(G')$. Furthermore, these points lie in some $\Tcal_p^+$.
    
    Denote by $\Bcal$ the right hull of $\Acal(G)\setminus \Ccal$ in  $\Acal(G')$.
    \begin{center}
    \begin{tikzpicture}
    \filldraw (0,0) circle (1pt) node[anchor=north] {$\Ccal$};
    \draw (-0.3,1) node[anchor=east] {$\mathcal{A}(G')$}; 
    \draw[dashed,ultra thick] (1,0) -- (4,0);
    \draw[decoration={brace,raise=5pt},decorate]
  (1,1) -- (5,1) ;
    \filldraw (0,1) circle (1pt);
    
    \draw (3,1.2) node[anchor=south] {$\mathcal{B}$};
    
    \draw (-0.3,0) node[anchor=east] {$\mathcal{A}(G)$};
    \draw[thick] (0.5,1) -- (5,1);
    \draw[dashed,ultra thick] (1,1) -- (4,1);
\end{tikzpicture}
\end{center}
    
    \begin{claim} $\Bcal=\Acal(H)\setminus\Ccal$. \end{claim} 
    \begin{proof}[Proof of Claim 2] By \Cref{cor:spines-of-cvx}, $\Acal(H)\setminus\Ccal$ is an end segment of $\Acal(G')$, and by definition of the right hull, so is $\Bcal$. It is clear that points in $\Bcal$ correspond to convex subgroups of $H$, that is, $\Bcal\subseteq\Acal(H)\setminus\Ccal$.
    
    For the other inclusion, take $a\in H_{>0}$. Since $H=\langle G\rangle$, there is $b\in G$ with $b>a$. In $\Acal(G')$, we have $\mathfrak t_p(b)\leqslant_{p,p}\mathfrak t_p(a)\leqslant_{p,p}\mathfrak s_p(a)$, that is, $\mathfrak t_p(a)$ and $\mathfrak s_p(a)$ are in $\Bcal$. As of $\mathfrak t_p^+(a)$, if it is not in $\Ccal$, then it defines a proper convex subgroup $G_{\mathfrak t_p^+(a)}$ of $H$ and again we can find $b\in G\setminus G_{\mathfrak t_p^+(a)}$, which gives $\mathfrak t_p(b)\geqslant\mathfrak t^+_p(a)$, and this concludes the claim.
\end{proof}

    $\Ccal$ is also a subset of $\Acal(H)$, since elements in $\Ccal$ arise from points in $G$ and $G\subseteq H$. Hence $\Acal(H)=\Ccal\cup\Bcal$.
    
    Now $(H,\Ccal + \Bcal)$, together with the restriction of $\mathfrak{t}_p,\mathfrak{s}_p, D^{[p^s]}_{p^r}, =_\bullet k_\bullet,\equiv_\bullet k_\bullet$ to $H$, is an $\Lcal_{syn}$-substructure of $G'$, and similarly $(G,\Acal(G))$ is an $\Lcal_{syn}$-substructure of $(H,\Ccal+\Bcal)$. By \Cref{lem:hulls!}, we have $\Acal(G)\setminus\Ccal\substruct\Bcal\substruct \Acal(G')\setminus \Ccal$. Thus, by \Cref{lem:sumofmo}, we get $$\Acal(G)=\Ccal+(\Acal(G)\setminus\Ccal)\substruct\Ccal+\Bcal=\Acal(H)\substruct \Ccal+(\Acal(G')\setminus\Ccal)=\Acal(G').$$

    By \Cref{fact:strongSchmitt}, we conclude $G\substruct H\substruct G'$.
\end{proof}

\begin{definition}
Let $G$ be an ordered abelian group. We say that $G$ is \emph{augmentable by infinites} (AI) if there exists a non-trivial $H$ such that $G\substruct H\oplus G$. Such an $H$ is called an \emph{infinite augment} of $G$.

Similarly, we say that $G$ is \emph{augmentable by infinitesimals} (Ai) if there exists a non-trivial $H$ such that $G\substruct G\oplus H$. Such an $H$ is called an \emph{infinitesimal augment} of $G$.
\end{definition}

Note that $G$ is augmentable by infinites if and only if $G$ admits a
proper elementary extension $G^*$ such that $G$ is a convex subgroup of $G^*$.
This equivalence is well known, and can in particular be derived from the
proof of the following theorem.

\begin{theorem}\label{thm:strongleftaugmentabilityOAG}
    All non-trivial ordered abelian groups are augmentable by infinites.
\end{theorem}

\begin{proof}
Let $G$ be a non-trivial ordered abelian group and let $G'$ be an elementary extension of $G$ realising a type at $+\infty$. Let $H$ be the convex hull of $G$ in $G'$. Then $G\substruct H$ by \Cref{lem:ConvxHull} and $G'/H\neq\{0\}$. Consider the exact sequence $0\rightarrow H\rightarrow G'\rightarrow G'/H\rightarrow 0$. Move if needed to an $\aleph_1$-saturated extension $(G^*,H^*)\superstruct (G',H)$. It is a well-known fact that now the corresponding exact sequence splits (see for example \cite{marker18}), that is, $G^*\simeq G^*/H^*\oplus H^*$. Since $H\substruct H^*$, then also $G\substruct H^*$. This implies that $A\oplus G\substruct A\oplus H^*$ for any ordered abelian group $A$, in particular, $G^*/H^*\oplus G\substruct G^*/H^*\oplus H^*=G^*$. Now since $G\substruct G^*$ and $G\subseteq G^*/H^*\oplus G$, we have that $G\substruct G^*/H^*\oplus G$.
\end{proof}

\begin{remark}\label{rmk:AI-is-elementary}
If $G\substruct H\oplus G$ and $H'\equiv H$, then $G\substruct H'\oplus G$, as $H\oplus G\equiv_G H'\oplus G$.
\end{remark}

\begin{corollary}
   Let $G$ be an ordered abelian group. Then for any prime $p$ TFAE:
   \begin{enumerate}
   \item Any infinite augment of $G$ is $p$-divisible,
   \item $T_p^+(G)$ has an initial point.
   \end{enumerate}

   Furthermore, TFAE:
   \begin{enumerate}
   \item $\mathbb Q$ is an infinite augment of $G$,
   \item Any divisible ordered abelian group is an infinite augment of $G$,
   \item For any $p$, $T_p^+(G)$ has an initial point.
   \end{enumerate}
   \label{cor_div}
\end{corollary}

\begin{proof}
An initial point $\alpha\in\Tcal_p^+(G)$ is such that $G_\alpha=G$, indeed, if not, then there is $a\in G\setminus G_\alpha$; but then $\mathfrak t_p^+(a)<\alpha$, which contradicts the definition of $\alpha$.

Assume that there is $a\in G$ such that $G_{\mathfrak t_p^+(a)}=G$. Let $Q$ be an infinite augment of $G$. Now since $G\substruct Q\oplus G$, $(0,a)\in Q\oplus G$ also defines the whole group, that is, $G_{\mathfrak t_p^+((0,a))}=Q\oplus G$. Assume towards a contradiction that $b\in Q$ is not $p$-divisible, then $(b,0)\in Q\oplus G$ is not $p$-divisible modulo $G$, thus $G\subseteq G_{\mathfrak s_p((b,0))}\subsetneq Q\oplus G$. But then $(0,a)\in G_{\mathfrak s_p((b,0))}$ and thus $G_{\mathfrak t_p^+((0,a))}\subseteq G_{\mathfrak s_p((b,0))}$, which means it can't equal the whole group $Q\oplus G$, a contradiction.

Conversely, assume that $G\substruct Q\oplus G$ for some non-trivial $p$-divisible ordered abelian group $Q$. Assume $(a,b)\in Q\oplus G$ is not $p$-divisible. Since $Q$ is $p$-divisible, $(a,b)$ is $p$-divisible modulo $G$, that is, $G_{\sfp((a,b))}\subseteq G$. This means that for any $a\neq0$, $G_{\tfpp((a,0))}$ must be equal to $Q\oplus G$, that is, $\Tpp$ has an initial point.

\ 

The second set of equivalences follows immediately from the first and from \Cref{rmk:AI-is-elementary} -- since the $\mathcal{L}_\textrm{oag}$-theory DOAG of divisible non-trivial ordered abelian groups is complete.
\end{proof}

\section{An application: definability of henselian valuations}\label{sec:def}


In this section, we give applications of Theorem \ref{thm:strongleftaugmentabilityOAG} to understand elementary embeddings of henselian valued fields, as well as (uniform) definability of henselian valuations. 
Recall that a field $k$ is called $t$-henselian if $k$ is $\mathcal{L}_\textrm{ring}$-elementarily
equivalent to a field admitting a non-trivial henselian valuation.
The key observation in what follows is that Theorem 
\ref{thm:strongleftaugmentabilityOAG} gives rise to an elementary embedding of any $t$-henselian
field of characteristic $0$ into some (non-trivial) power series field over it:

\begin{proposition}
    Let $k$ be a $t$-henselian field of characteristic $0$. Then 
    $k \substruct k((H))$ for some non-trivial ordered abelian group $H$. \label{power}
\end{proposition}
\begin{proof}
    Let $k$ be a $t$-henselian field of characteristic $0$. 
    We first argue that
    we may replace $k$ with any $\mathcal{L}_\textrm{ring}$-elementary extension $k^*$: 
    Assume that $k\substruct k^*\substruct k^*((H))$ holds for some non-trivial ordered abelian group $H$.
    We obtain elementary embeddings
    $k \substruct k^*((H))$ and 
    $k\subseteq  k((H)) \substruct k^*((H))$, where the latter embedding is elementary 
    by the Ax-Kochen/Ershov Theorem in equicharacteristic $0$ \cite[Theorem 6.17]{Hils-mtvf}. These embeddings
    imply that $k\substruct k((H))$ also holds.
    
    In particular, by the paragraph above, 
    we may assume that $k$ is $\aleph_1$-saturated in $\mathcal{L}_\textrm{ring}$.
    By \cite[Lemma 3.3]{PZ78}, $k$ admits a non-trivial henselian valuation $v$. 
    Passing to an $\aleph_1$-saturated $\mathcal{L}_\textrm{val}$-elementary extension $(k^*,v^*)$ of $(k,v)$, we may further assume that $(k,v)$ is $\aleph_1$-saturated as a valued field.
    We next argue that $k$ admits a non-trivial henselian valuation $w$
    of residue characteristic $0$. 
    Note that $v$ has a non-trivial coarsening $w$ of residue characteristic $0$ if and only if 
    $vk\supsetneq \braket{v(p)}_{convex}$ holds. That this inclusion is proper is ensured by $\aleph_1$-saturation of $(k,v)$. Thus, replacing $v$ by $w$, we may assume $\mathrm{char}(kv)=0$.

    We now argue that we have $(k,v)\substruct (kv((vk)),u)$, where $u$ denotes the power series
    valuation on $kv((vk))$. 
    This argument is also in \cite{vdD14}, and we include
    it for the convenience of the reader.
    Since  $(k,v)$ be a $\aleph_1$-saturated, there is a section of the valuation $s:vk \rightarrow k$. Since $(k,v)$ is of equicharacteristic $0$, there is a lift $l:kv\rightarrow k$ of the residue field. In particular, it follows that $kv(vk)$ embeds in $k$. The extension is immediate and therefore, $k$ and $kv(vk)$ admit $kv((vk))$ for unique (up to isomorphism) maximal extension. Then we have  $(k,v) \substruct (kv((vk)),u)$ by Ax-Kochen/Ershov.

    Let now $G=vk$. By Theorem \ref{thm:strongleftaugmentabilityOAG}, we have
    $G \substruct H \oplus G$ for some non-trivial ordered abelian group $H$.
    Thus, we obtain -- once again by the Ax-Kochen/Ershov Theorem in equicharacteristic
    $0$ --
    $$k \substruct kv((G)) \substruct kv((G))((H))$$
    and $$k((H)) \substruct kv((G))((H)).$$
    This induces an elementary embedding $k \substruct k((H))$.
\end{proof}

\begin{remark}
The proposition does not hold for fields of positive characteristic in general.
Note that if $k$ is perfect, then $k \substruct k((H))$ implies that $H$ is $p$-divisible and hence the power series valuation on $k((H))$ is tame.
However, there are examples of perfect henselian fields which have no elementary extension
that admits a non-trivial tame
valuation.
One such example is
$k=(\mathbb{F}_p((t))^\mathcal{U})^\mathrm{perf}$, the perfect hull of any nonprincipal
ultrapower of the power series field $\mathbb{F}_p((t))$.
\label{rem_Jonas}
\end{remark}

We now link the proposition above to $\mathcal{L}_\textrm{ring}$-definability of henselian
valuations with a given residue field. 
Recall that a valuation $v$ on a field $K$ is \emph{definable} if there is
an $\mathcal{L}_\textrm{ring}$-formula $\psi(x)$, possibly using parameters from $K$, such that
$\psi(K)=\mathcal{O}_v$ holds. If $v$ is definable via a formula $\psi$ that requires 
no parameters from $K$, we say it is \emph{$\emptyset$-definable}.

A key fact to establish definability abstractly is the following:
\begin{fact}[Beth definability theorem, {\cite[Thm.~6.6.4]{Hodges}}]
    Let $\Lcal$, $\Lcal'$ be first-order languages with $\Lcal\subseteq\Lcal'$. Let $T$ be a theory in $\Lcal'$ and $\phi(\bar{x})$ a formula of $\Lcal'$. Then the following are equivalent:
    \begin{itemize}
        \item If $A$ and $B$ are models of $T$ such that $\restriction{A}{\Lcal}=\restriction{B}{\Lcal}$, then $A\models\phi(\bar{a})$ if and only if $B\models\phi(\bar{a})$, for all tuples $\bar{a}\in A$;
        \item $\phi(\bar{x})$ is equivalent modulo $T$ to a formula $\psi(\bar{x})$ of $\Lcal$.
    \end{itemize} \label{fact:beth}
\end{fact}


\begin{corollary}
     Let $(K,v)$ be a henselian valued field, with $Kv$ of characteristic $0$.  
    Then the following are equivalent:
    \begin{enumerate}
        \item $v$ is \emph{not} definable (over parameters) in the language of rings,
        \item there is an elementary extension $(K^*,u)$ of $(K,v)$ and a valuation $w\neq u$
        on $K^*$ such that $(K^*,u) \equiv_K (K^*,w)$. 
    \end{enumerate} \label{comparable}
\end{corollary}
\begin{proof} The equivalence 
follows immediately from Beth definability theorem, applied with the language $\lL=\lL_\textrm{ring}(K)$ (with constant symbol for each element in $K$), 
$\lL'=\lL_\textrm{ring}(K)\cup\{\mathcal{O}\}$ 
 (with an additional predicate for the valuation ring),
 the theory $T=\Th_{\lL'}(K,v)$ and where $\phi(x)$ is the $\lL_\textrm{val}$-formula $x\in \mathcal{O}$ 
 defining the valuation ring in $T$. 
\end{proof}

\begin{theorem} \label{def_fieldwise}
    Let $k$ be a field of characteristic $0$. The following are equivalent:
    \begin{enumerate}
        \item $k$ is not t-henselian, i.e., it is not elementary equivalent to a field $k'$ that admits a non-trivial henselian valuation,
        \item for every henselian valued field $(K,v)$ with residue field $k$, the valuation $v$ is definable (possibly using parameters),
        \item for every henselian valued field 
        $(K,v)$ with residue field $k$, the valuation $v$ is 
        $\emptyset$-definable,
        \item there is a parameter-free $\mathcal{L}_\textrm{ring}$-formula $\psi(x)$ which defines the valuation ring $v$ of in any henselian valued field $(K,v)$ with residue field a model of $\mathrm{Th}_{\mathcal{L}_\textrm{ring}}(k)$. In other words,
        all henselian valuations with residue field elementarily equivalent to $k$ are
        uniformly $\emptyset$-definable in $\mathcal{L}_\textrm{ring}$.
    \end{enumerate}
\end{theorem}
\begin{proof} The implications (4) $\implies$ (3) and (3) $\implies$ (2) are trivial.

    We show (2) $\implies$ (1) via contraposition.
    Assume that $k$ is $t$-henselian. By Proposition \ref{power},
    we have $k \substruct k((H))$ for 
    some non-trivial ordered abelian group $H$. 
    Notice that we have the following chain of elementary embeddings :
    \[k \substruct k((H_0)) \substruct   k((H_{-1} + H_0 + H_{1})) \substruct \cdots \substruct k((\sum_{-n<i<n} H_i)) \substruct \cdots\]
    where $H_i$ are all copies of $H$ and with the obvious embeddings. By Tarski's chain lemma \cite[Theorem~2.1.4]{TZ12}, $k$ is an elementary substructure of $k':=k((\sum_{i \in 
    \Zbb}H))$. We show that there is a henselian valued field with residue $k'$ whose valuation is not definable in $\lL_{ring}$. The same property will hold for $k$ instead of $k'$ by the $\preceq$-version of the  
    Ax-Kochen/Ershov theorem (\cite[Theorem 6.17]{Hils-mtvf}).
    Now consider  $\Gamma:=\bigoplus_{i \in \Zbb} H$
   and the Hahn field 
   $K=k'((\Gamma))$, together with the power series
   valuation $v$ with value group $\Gamma$. Consider the field
   $K^* = k'((\Gamma'))((\Gamma))$ with $\Gamma'$ is another copy of $\Gamma$.  
   Note that $\Gamma\oplus\Gamma'$ is an elementary extension of $\Gamma$
   (this holds, e.g., by playing Ehrenfeucht–Fraïssé games). The field $K^*$
   admits the distinct henselian valuations 
   \begin{itemize}
       \item $u$ with residue field $k'$ and value group $\Gamma \oplus \Gamma'$ and
       \item $w$ with residue field $k'((\Gamma'))$ and value group $\Gamma$.
   \end{itemize}
   Then we have by applying Ax-Kochen/Ershov \cite[Theorem 6.17]{Hils-mtvf}
   once again
   $$(K,v)\substruct (K^*,u) \textrm{ and }(K,v)\substruct (K^*,w)$$
   for the natural (field) embedding of $K$ in $K^*$.
    By Corollary \ref{comparable}, this shows that $v$ is not $\mathcal{L}_\textrm{ring}$-definable on $K$.

    The implication (1) $\implies$ (4) is essentially \cite[Proposition 5.5]{FJ15}. 
    Assume $k$ is not $t$-henselian, in particular $k$ is not separably closed.
    Consider the $\mathcal{L}_\textrm{val}$-theory $T$ which stipulates that any model $(K,v)$ is henselian with residue field $Kv \equiv k$ (in $\mathcal{L}_\textrm{ring}$) and the 
    $\mathcal{L}_\textrm{val}$-formula $\phi(x)$ which asserts that $ x \in \mathcal{O}_v$.
    Take $(K_1,v_1), (K_2,v_2) \models T$ that have the same $\mathcal{L}_\textrm{ring}$-reduct, i.e., in particular $K_1=K_2$. Then both $v_1$ and $v_2$ are both henselian valuations with
    non-separably closed residue field, thus they are comparable, cf.~\cite[Theorem 4.4.2]{EP05}. 
    Assume for a contradiction that we have $\mathcal{O}_{v_1} \subsetneq \mathcal{O}_{v_2}$, then $v_1$ induces a
    non-trivial henselian
    valuation on $Kv_2$. But then $k$ is $t$-henselian, in contradiction to our assumption. Thus,
    we have $v_1=v_2$, and so Fact \ref{fact:beth} implies that $v$ is uniformly $\emptyset$-definable in
    all models of $T$.
\end{proof}

\begin{remark} The crucial ingredient in
Theorem \ref{def_fieldwise} is Proposition \ref{power}, which fails in positive characteristic
in general
(cf.~Remark \ref{rem_Jonas}).
However, if $k$ is a perfect field, amending the proof of Theorem \ref{def_fieldwise} 
(using \cite[Theorem 7.1]{Kuh16} instead of \cite[Theorem~6.17]{Hils-mtvf}) shows we have the following chain of implications:

\begin{center}$k \substruct k((\Gamma))$ for some non-trivial
ordered abelian group $\Gamma$\\
$\Downarrow$\\ 
there is an henselian valued field $(K,v)$ with $Kv=k$ such that $v$ is
not definable\\
$\Downarrow$\\ 
there is an henselian valued field $(K,v)$ with $Kv=k$ such that $v$ is
not $\emptyset$-definable\\
$\Downarrow$ \\ 
henselian valuations with residue field elementarily equivalent to
$k$ are not uniformly $\emptyset$-definable
in $\mathcal{L}_\textrm{ring}$\\
$\Downarrow$\\
$k$ is $t$-henselian.\\
\end{center}
We do not know whether Theorem \ref{def_fieldwise} holds or fails
for perfect fields $k$ in general. Since Proposition \ref{power} may fail, 
our proof method certainly does not apply.

Note that Theorem \ref{def_fieldwise} may fail for imperfect fields: 
If $k$ is imperfect and admits no Galois extensions of degree divisible by $p$, then any henselian valuation with residue field $k$ is definable by \cite[Proposition 3.6]{Jah24}.
Note that such examples exist abundantly, e.g.,
$\mathbb{F}_p(t)^\mathrm{sep}$ is t-henselian (and in fact even henselian) but admits no proper
Galois extensions of any degree.
\end{remark}

Theorem \ref{def_fieldwise} partially answers a question raised by Krapp, Kuhlmann, and Link (\cite[\S 5.4]{KKL}). Specifically, they inquire about identifying an optimal class $\Ccal$ of fields such that for any
henselian valued field with residue field in $\Ccal$, the valuation is definable in the language of rings. Moreover,
they ask the analogous question about identifying the subclass $\Ccal_0$ of $\Ccal$ such that
such that for any henselian valued field with residue field in $\Ccal$, the valuation is 
$\emptyset$-definable.
Theorem \ref{def_fieldwise} shows that restricting to the case of characteristic $0$,
the class of fields that are not $t$-henselian coincides
exactly with both the classes $\Ccal$ and $\Ccal_0$. 


\renewcommand{\appendixname}{}
\appendix
\section{Iksrat's Lemma and another proof of augmentability}\label{sec:app}

Before coming up with our (much simpler!) argument to classify augmentable multi-orders, our strategy involved considering the intersection of a descending elementary chain of models and claiming it is a model. Surprisingly, to the extent of our knowledge, this was not a well-known fact. We decided to present this result here, as an appendix, despite the fact that it is currently a lemma with no implications, except another (more complicated) proof of \Cref{lem:aug-mo}.

\begin{theorem}[Iksrat's lemma]\label{iksrat} Let $(I,<)$ be a linear order, and let $M$, $N$, $(A_i)_{i\in I}$ be structures in a first-order language $\Lcal$. Assume that $M\substruct A_j\substruct A_i\substruct N$ for any $i<j$. Let $\kappa=\cf(I)$ and assume furthermore that $N$ is $\kappa^+$-saturated. 
Let $A=\bigcap_{i\in I}A_i$. Assume finally that for each $i\in I$ and for any $a_i\in A_i$, there is an $A$-definable set $B_i$ (in $N$) such that $a_i\in B_i\subseteq A_i$.

Then $M\substruct A\substruct N$.
\end{theorem}

Note that the existence of $A$-definable $A_i$-neighbourhoods of any element $a_i\in A_i$ holds in particular if $A_i$ itself is $M$-definable. We also note that other definability conditions on $A_i$ give similar results, for example, assuming that there exists a $N$-definable set $B_i$ such that $(A_j\setminus A)\subseteq(B_i\setminus A)\subseteq(A_i\setminus A)$.

\begin{proof}
 If $\kappa=1$, $I$ has a final element $i$ and $A=A_{i}$, so the lemma obviously holds.
 
 If $\kappa$ is a limit ordinal, we take a sequence $(i_j)_{j\in\kappa}$ cofinal in $I$. This means $A=\bigcap_{j\leqslant\kappa} A_{i_j}$. In order to keep our notation readable, we will only work with this sequence and identify $A_j$ with $A_{i_j}$. We prove by induction on $\Lcal$-formulas $\phi$ the following result: for any $\overline a\in A$, $A\satisf\phi(\overline a)$ iff $N\satisf\phi(\overline a)$. This will give $A\substruct N$, from which follows $M\substruct A$ since $M\substruct N$ and $M\subseteq A$.
 
 If $\phi$ is atomic, the result holds. If $\phi$ is a negation or a conjunction, the result holds by induction. Now let $\phi(\overline x)=\exists y\psi(\overline x,y)$ and assume for induction that the result holds for $\psi$. Let $\overline a\in A$.
 
 If $A\satisf\phi(\overline a)$, then there is $a'\in A$ such that $A\satisf\psi(\overline a,a')$, thus by induction $N\satisf\psi(\overline a,a')$ and therefore $N\satisf\phi(\overline m)$.
 
 Conversely, if $N\satisf\phi(\overline a)$, then, since $N\superstruct A_i$ holds for all $i<\kappa$, we have $A_i\satisf\phi(\overline a)$. We will exhibit, by induction on $i$, $A$-definable sets $B_i\subseteq A_i$ such that for all $i<\kappa$, $B_{i+1}\subseteq B_i$ and there is $a_i\in\bigcap_{j\leqslant i}B_i$ with $A_i\satisf\psi(\overline a,a_i)$. 

 We know $A_0\satisf\phi(\overline a)$, thus there is $a_0\in A_0$ such that $A_0\satisf\psi(\overline a,a_0)$, and by assumption of the lemma there is an $A$-definable set $B_0\subseteq A_0$ such that $a_0\in B_0$.
 
 If $i=j+1$, assuming $(B_k)_{k<i}$ has been constructed, then $A_j\satisf\exists y(\psi(\overline a,y)\et y\in B_j)$. Since $A_i\substruct A_j$ and $B_j$ is $A$-definable, $A_i\satisf\exists y(\psi(\overline a,y)\et y\in B_j)$. Let $a_i$ witness this formula and let $B'_i\subseteq A_i$ be an $A$-definable set such that $a_i\in B'_i$. Take $B_i=B_i'\cap B_j$, it contains $a_i$, is $A$-definable and is contained in $B_j$, as wanted.
 
 If $i$ is a limit ordinal, assuming $(B_j)_{j<i}$ has been constructed, for each $j<i$ since $A_j\satisf\exists\psi(\overline a,y)\et y\in B_j$ and $A_i\substruct A_j$, there is $a_{i,j}\in B_j\cap A_i$ such that $A_i\satisf\psi(\overline a,a_{i,j})$. Let $B'_{i,j}\subset A_i$ be an $A$-definable set containing $a_{i,j}$ and let $B_{i,j}=B_j\cap B'_{i,j}$. Consider the type $p(y)=\sset{\psi(\overline a,y)\et y\in B_{i,j}}{j<i}$, its restriction $p_{<k}(y)=\sset{\psi(\overline a,y)\et y\in B_{i,j}}{j<k}$ is realised by $a_{i,k}$, thus it is realised by some $a_i$ in $N$ by saturation. But now $a_i\in\bigcap_{j<i} B_{i,j}\subseteq A_i$, so $N\satisf\psi(\overline a,a_i)$, and since $A_i\substruct N$ also $A_i\satisf\psi(\overline a,a_i)$. Now take $B_i\subseteq A_i$ $A$-definable containing $a_i$. Since $a_i\in\bigcap_{j<i} B_{i,j}$, in particular, $a_i\in\bigcap_{j\leqslant i}B_j$.
 
 Finally, consider the type $p(y)=\sset{\psi(\overline a,y)\et y\in B_i}{i<\kappa}$. By construction of $(B_i)_{i\in I}$, it is finitely satisfied, thus it is realised by some $a'\in N$. Now $a'\in\bigcap_{i<\kappa} B_i\subseteq\bigcap_{i<\kappa} A_i=A$, so $N\satisf\psi(\overline a,a')$ and by induction (on formulas) also $A\satisf\psi(\overline a,a')$, that is, $A\satisf\phi(\overline a)$.
\end{proof}

\begin{remark}
Iksrat's lemma is not quite simply backward Tarski's chain lemma, since the saturation and definability assumptions are important:
\begin{itemize}
\item in the first-order language of pure orders, Let $N=\mathbb R$, $M=\mathbb R_{<0}$ and $A_i=\mathbb R_{<\tfrac{1}{i}}$ with $i\in\mathbb N$. We have $M\substruct A_i\substruct N$, each $A_i$ is definable, but $N$ is not $\aleph_1$-saturated; and $A=\bigcap_{i\in \mathbb N}A_i =\mathbb R_{\leqslant 0}$ is not a model since it has a max.
\item similarly, take $N=\mathbb R^*\superstruct\mathbb R$\ \ $\aleph_1$-saturated, but take $M=\mathbb R_{<0}$ and $A_i=\mathbb R_{<\tfrac{1}{i}}$ as above; they still only contain standard reals and thus are very much not definable in $N$ themselves, nor is any subset of them. Again $M\substruct A_i\substruct N$ but $A=\mathbb R_{\leqslant 0}$ has a max and is not a model.
\end{itemize}
\end{remark}

\begin{corollary}[\Cref{lem:aug-mo}]
Let $A$ be a coloured $P$-multi-order. $A$ is augmentable on the right iff for some $p\in P$, $A_p$ is unbounded on the right. \label{cor_iksrat}
\end{corollary}

\begin{proof}
Let $A\substruct A^*$ be $|A|^+$-saturated. For each $q\in P$ and each $a\in A_q$, Let $L_a=\sset{b\in A_r}{b\leqslant_{r,q}a}_{r\in P}$ and $R_a=\sset{b\in A_r}{b\leqslant_{q,r}a\wedge\neg b\leqslant_{r,q}a}_{r\in P}$. Define similarly $L_a^*$ and $R_a^*$ as subsets of $A^*$. Finally, let $C_a=L_a+R_a^*$ and $C=\bigcap_{a\in A} C_a$.

We claim that $A\substruct C_a\substruct A^*$ for all $a\in A$ and that $C_a\substruct C_b$ if $b\leqslant_{q,r} a$. This is obvious from \Cref{lem:sumofmo}.

Now fix $a\in A$ and let $b\in C_a$. If $b\in A$, then $\sset{b}$ is $A$-definable. If $b\notin A$, then $b\in R^*_a$. For some $q\in P$, $b\in A^*_q$, and $(R^*_a)_q=\sset{x\in A^*_q}{a\leqslant_{r,q}x\wedge\neg x\leqslant{q,r}a}$ is $a$-definable. In both cases, we found an $A$-definable neighbourhood of $b$ included in $C_a$.

Hence we can apply Iksrat's lemma (\Cref{iksrat}) with $M=A$, $N=A^*$, and $I=A$; we therefore have $A\substruct C\substruct A^*$. Now, $C=A+\bigcap_{a\in A} R^*_a$, which means $A\substruct A+\bigcap_{a\in A} R^*_a$; of course $\bigcap R^*_a\neq\emptyset$ since it contains a realisation of a type at infinity.
\end{proof}

\bibliographystyle{abbrv}
\bibliography{bibtex}

\end{document}

%% file: preamble.tex
\usepackage[T1]{fontenc} 
\usepackage[utf8]{inputenc}
\usepackage{marvosym}
\usepackage{amsmath}
\usepackage{scalerel}
\usepackage{amsthm}
\usepackage{amsfonts}
\usepackage{amssymb}
\usepackage{comment}
\usepackage{graphicx}
\usepackage{slashed} 
\usepackage{braket}
\usepackage{hyperref}
\usepackage[nameinlink]{cleveref}
\usepackage{tikz-cd}
\usepackage{adjustbox}
\usepackage{euscript}
\usepackage{thmtools}
\usepackage{thm-restate}
\usepackage{xparse}
\usepackage{mathtools}
\usetikzlibrary{decorations.pathreplacing,angles,quotes}

\usepackage[inline]{enumitem}

\hypersetup{
  colorlinks   = true, 
  urlcolor     = {blue!50!black}, 
  linkcolor    = {blue!50!black}, 
  citecolor   = {red!50!black} 
}

\def\restriction#1#2{\mathchoice
              {\setbox1\hbox{${\displaystyle #1}_{\scriptstyle #2}$}
              \restrictionaux{#1}{#2}}
              {\setbox1\hbox{${\textstyle #1}_{\scriptstyle #2}$}
              \restrictionaux{#1}{#2}}
              {\setbox1\hbox{${\scriptstyle #1}_{\scriptscriptstyle #2}$}
              \restrictionaux{#1}{#2}}
              {\setbox1\hbox{${\scriptscriptstyle #1}_{\scriptscriptstyle #2}$}
              \restrictionaux{#1}{#2}}}
\def\restrictionaux#1#2{{#1\,\smash{\vrule height .8\ht1 depth .85\dp1}}_{\,#2}}

\DeclareMathOperator{\discr}{discr}

\DeclareMathOperator{\Th}{Th}

\DeclareMathOperator{\cf}{cf}

\makeatletter
\DeclareRobustCommand\bigop[1]{%
  \mathop{\vphantom{\sum}\mathpalette\bigop@{#1}}\slimits@
}
\newcommand{\bigop@}[2]{%
  \vcenter{%
    \sbox\z@{$#1\sum$}%
    \hbox{\resizebox{\ifx#1\displaystyle.9\fi\dimexpr\ht\z@+\dp\z@}{!}{$\m@th#2$}}%
  }%
}
\makeatother

\newcommand{\Acal}{\ensuremath{\mathcal{A}}}
\newcommand{\Bcal}{\ensuremath{\mathcal{B}}}
\newcommand{\Ccal}{\ensuremath{\mathcal{C}}}

\newcommand{\Lcal}{\ensuremath{\mathcal{L}}}

\newcommand{\Scal}{\ensuremath{\mathcal{S}}}
\newcommand{\Tcal}{\ensuremath{\mathcal{T}}}

\newcommand{\Nbb}{\ensuremath{\mathbb{N}}}

\newcommand{\Zbb}{\ensuremath{\mathbb{Z}}}

\newcommand\cla\EuScript
\let\lang\CMcal

\newcommand{\lL}{\lang{L}}

\usepackage{ulem}
\usepackage{todonotes}
\setuptodonotes{noinline}

\colorlet{PierreColor}{cyan}

\colorlet{BlaiseColor}{violet!80!white}

\colorlet{AnnaColor}{green!80!yellow}

\colorlet{FranziColor}{yellow}

\NewDocumentCommand{\sset}{mg}{\left\{#1\IfNoValueF{#2}{\mid#2}\right\}}

\theoremstyle{plain}
\newtheorem{theorem}{Theorem}[section]

\newtheorem{lemma}[theorem]{Lemma}
\newtheorem{proposition}[theorem]{Proposition}
\newtheorem{fact}[theorem]{Fact}
\newtheorem{corollary}[theorem]{Corollary}

\theoremstyle{definition}
\newtheorem{definition}[theorem]{Definition}

\newtheorem*{notation*}{Notation}

\newtheorem*{question*}{Question}
\newtheorem{example}[theorem]{Example}

\theoremstyle{remark}

\newtheorem{remark}[theorem]{Remark}
\newtheorem*{remark*}{Remark}

\newtheorem{claim}{Claim}
\newtheorem*{claim*}{Claim}
\newtheorem{poc}{Proof of Claim}

\definecolor{black}{rgb}{0,0,0}

\colorlet{savedColor}{.}

\crefname{subsection}{subsection}{subsections}
\crefname{subsubsection}{subsubsection}{subsubsections}

\let\et\wedge

\let\satisf\vDash
\let\substruct\preccurlyeq
\let\superstruct\succcurlyeq